\documentclass[11pt,letterpaper]{article}

\usepackage{url}
\usepackage{tikz}
\usepackage{fullpage}

\usepackage[subtle,mathspacing=normal]{savetrees}
\usepackage{amsmath}
\usepackage{algorithm}
\usepackage[noend]{algpseudocode}
\usepackage{amssymb}
\usepackage{amsthm}
\usepackage{color}
\usepackage{array}
\usepackage{xy}
\usepackage{setspace}
\usepackage{varwidth}
\usepackage{multicol}
\usepackage{algorithm}
\usepackage{hyperref}
\usepackage{thm-restate}

\newtheorem{thm}{Theorem}[section]
\newtheorem{theorem}{Theorem}[section]
\newtheorem{definition}[thm]{Definition}
\newtheorem{lemma}[thm]{Lemma}

\newtheorem{claim}[thm]{Claim}

\theoremstyle{remark}

\newtheorem{problem}[thm]{Problem}

\renewcommand{\paragraph}[1]{\vspace{.2 cm} \noindent \textbf{#1}}

\usepackage{authblk}

\def\sub{\subseteq}

\def\vv{{\bf v}}
\def\vec#1{{\bf v}}
\def\eps{\varepsilon}
\def\path#1{\langle{#1}\rangle}
\def\edge#1#2{\langle{#1,#2}\rangle}

\begin{document}

\title{Variants of Baranyai's Theorem with Additional Conditions}
\author{Zoe Xi}
\date{\today}
\maketitle

\begin{abstract}
A classical theorem of Baranyai states that, given integers $2\leq k <
n$ such that $k$ divides $n$, one can find a family of ${n-1\choose k-1}$
partitions of $[n]$ into $k$-element subsets such that every subset
appears in exactly one partition. In this paper, we build on recent work by Katona and Katona in studying partial partitions, or parpartitions, of $[n]$ that consist of $k$-element sets not
overlapping significantly.  More precisely, two parpartitions $P_1$
and $P_2$ are considered $(\alpha,\beta)$-close for $\alpha,\beta\in
(0,1)$ if there exist subsets $A_1\neq B_1\in P_1$ and $A_2\neq
B_2\in P_2$ such that $|A_1\cap A_2| > \alpha{k}$ and $|B_1\cap B_2| >
\beta{k}$.


We establish that, given integers $k$, $\ell$, and $n$ satisfying
$k^2\ell\leq n$ and $\alpha, \beta\in (0, 1)$ satisfying
$\alpha+\beta\geq{(k+2)/k}$, one can find $\lfloor {n\choose
  k}/\ell\rfloor$ $(k, \ell)$-parpartitions of $[n]$ such that no two
distinct $(k, \ell)$-parpartitions are $(\alpha,\beta)$-close; this
result improves the condition $k=O(1)$ and $\ell=o(\sqrt{n})$ in a
corresponding result by Katona and Katona for $\alpha = \beta = 1/2$.
We also prove that, given integers $k$, $\ell$, and $n$ satisfying
$k=O(1)$ and $\ell=o(\sqrt{n})$, there is a cyclic ordering of the
$k$-element subsets of $[n]$ for any chosen $\alpha+\beta\geq{1}$ such
that any $\ell$ consecutive $k$-element subsets in the ordering form a
$(k, \ell)$-parpartition of $[n]$, which we refer to as a consecutive
$(k, \ell)$-parpartition (according to the ordering), and any two of
these disjoint consecutive $(k, \ell)$-parpartitions are not
$(\alpha,\beta)$-close.
\end{abstract}


\section{Introduction}
\label{section:introduction}
We begin by stating the following classical theorem of Baranyai
\cite{baranyai1974factrization}.

\begin{theorem}[Baranyai~1973]
\label{thm:baranyai}
Let $2\leq k < n$ be integers such that $k$ divides $n$. There exist
${n-1\choose k-1}$ partitions of $[n]$ into $k$-element subsets such
that no subset appears in two of these partitions.
\end{theorem}  

Baranyai's theorem has applications to areas including database
theory~\cite{demetrovics1998design}, coding
theory~\cite{zhang2014coding}, and information
theory~\cite{tamm1996applications}. In an application, one would often
like to say something further about the relationship between
partitions (consisting of $k$-element subsets), which Baranyai's
theorem places no restriction on~\cite{katona2005constructions}. A
natural research direction is to consider Baranyai's theorem (or some
weaker version) under the additional restriction requiring that no two
distinct partitions $P_1$ and $P_2$ should be ``close.''  For
instance, we say that two partitions $P_1$ and $P_2$ (consisting of
$k$-element subsets) are \emph{$(\alpha,\beta)$-close} for
$\alpha,\beta\in(0,1)$ if there exist subsets $A_1\neq B_1\in P_1$ and
$A_2\neq B_2\in P_2$ such that
$|A_1\cap A_2| > \alpha{k}$ and $|B_1\cap B_2| > \beta{k}$. We write that $P_1$ and
$P_2$ are \emph{close} if they are $(\alpha,\beta)$-close
for some $\alpha+\beta\geq{1}$.  According to this definition of
closeness, if there is a pair of subsets $(A_1, B_1)$ in $P_1\times
P_2$ with large intersection, then all other disjoint pairs in
$P_1\times P_2$ must have small intersection in order for $P_1$ and
$P_2$ to be not close.

A recent paper by Katona and
Katona~\cite{katona2023towards} studies this restriction for
the specific case of $\alpha=\beta=1/2$ and for the concept of a
\emph{$(k, \ell)$-partial partition}, or \emph{parpartition}, of $[n]$, which
is a family of $\ell$ disjoint $k$-element subsets of $[n]$ whose union might be a strict subset of $[n]$ (where
$k\ell\leq{n}$ holds). The motivation for replacing partitions with
parpartitions is that, as Katona and Katona write, this problem for
partitions seems very difficult, but becomes tenable for the weaker
concept of parpartitions. 

\paragraph{Our results.}
Our first result is a variant of Baranyai's theorem for parpartitions that are not too close. 

\begin{theorem}\label{thm:modified-katona}
Let $n$ be a positive integer. Fix a constant $k$, and let $\ell$ be such that $k^2\ell\leq
n/3$ and $m={n\choose k}$.
Let $\alpha$, $\beta\in (0, 1)$ such that $\alpha+\beta\geq (k+2)/k$. One can find $\lfloor{m}/{\ell}\rfloor$ $(k, \ell)$-parpartitions
of $[n]$ such that 
\begin{itemize}

\item[(1)] no $k$-element subset appears in two
parpartitions and

\item[(2)] no two parpartitions $P_1$ and $P_2$ contain sets $A_1, B_1\in P_1$ and $A_2, B_2\in P_2$
satisfying $|A_1\cap A_2|>\alpha k$ and $|B_1\cap B_2|>\beta k$.
    
\end{itemize}
\end{theorem}

This result should be compared with a corresponding theorem by Katona and Katona \cite{katona2023towards} that we state as follows.

\begin{theorem}[Theorem 2.1, Katona and Katona 2023]
\label{thm:katona-katona-1}
Let $k$, $\ell$, and $n$ be
positive integers such that $k$ is some fixed constant and $\ell=o(\sqrt{n})$. Let $m={n\choose
  k}$. If $n$ is large enough, one can find
$\lfloor{m}/\ell\rfloor$ $(k, \ell)$-parpartitions of $[n]$ into
$k$-element subsets such that no subset appears in two parpartitions and
no two parpartitions are $(1/2,1/2)$-close.
\end{theorem}
Note that the condition that $k$ is some fixed constant and $\ell=o(\sqrt{n})$ here is based on the proof of the theorem given in \cite{katona2023towards}.

Define the \emph{size} of a $(k, \ell)$-parpartition of $[n]$ to be
$k\ell$. Theorem \ref{thm:modified-katona} increases the size of
a parpartition in Theorem \ref{thm:katona-katona-1} from $o(\sqrt{n})$
to $\Theta(n)$ 
with the small tradeoff that the threshold in the definition
of $(\alpha, \beta)$-close is strengthened to $\alpha+\beta\geq
(k+2)/k$ from $\alpha+\beta\geq 1$. Note that $\alpha$ and $\beta$ are fixed to be $1/2$ in the result by Katona and Katona. 



To prove their result, Katona and Katona take an indirect approach, solving a graph problem that
we paraphrase in Problem \ref{problem:graph-reduction}. Let $G_1 = (V, E_1)$ and $G_2 = (V, E_2)$ be
two simple graphs on the same vertex set $V$.
%
%
Define an \emph{alternating-$(\ell, 2, \ell,
  2)$-bag} to consist of two vertex-disjoint copies $A$ and $B$ of
$K_\ell$, the clique of size $\ell$, in $G_1$ and two vertex-disjoint
edges in $G_2$, each connecting a vertex in $A$ to a vertex in $B$. As
there are two graphs $G_1$ and $G_2$ involved in this definition, we
refer to this notion of an $(\ell, 2, \ell, 2)$-bag as a two-graph version.
We will give a three-graph version of an $(\ell, 2, \ell, 2)$-bag later.

\begin{problem}
\label{problem:graph-reduction}
Let $G_1 = (V, E_1)$ and $G_2 = (V, E_2)$ be two simple graphs on the
same vertex set $V$ such that their edge sets $E_1$ and $E_2$ are
disjoint. Let $\ell\geq 2$ be an integer. Let $\delta_1$
(resp.~$\Delta_2$) denote the minimum (resp.~maximum) degree of a
vertex in $G_1$ (resp.~$G_2$). Impose some appropriate conditions on
$\delta_1$ and $\Delta_2$ such that, if $m$ is large enough, then
there are $\lfloor m/\ell\rfloor$ vertex-disjoint copies of $K_\ell$
in $G_1$ such that no two of these copies span an alternating-$(\ell,
2, \ell, 2)$-bag in $G_1$ and $G_2$; these conditions can be
discharged later when the result is applied.
\end{problem}  

Our next result answers a problem posed by Katona and Katona, which is
essentially Problem~\ref{problem:graph-reduction} extended in a
direction involving certain powers of Hamiltonian cycles. The nature
of the problem is somewhat technical. Given a Hamiltonian cycle $H$
(that is, a cycle that passes through every vertex in a graph exactly
once), the $p$th power $H^p$ of $H$ is obtained by adding an edge
between every pair of vertices connected by a path of length at most
$p$ in $H$. One can readily observe that $H^p$ contains $m$ copies of
$K_{p+1}$ for every natural number $p$ less than the size of $H$.

Note that, in this result, {alternating-$(\ell, 2, \ell, 2)$-bag} is a three-graph
version, which consists of two vertex-disjoint copies $A$ and $B$ of
$K_\ell$ in $G_1$ and two vertex-disjoint edges, one in ${G_2}$ and
the other in ${G_3}$, each connecting a vertex in $A$ to a vertex in
$B$. For the sake of a clearer presentation, we assume in the rest of
the paper that the colors red, green, and blue are assigned to the edges in
$G_1$, $G_2$, and $G_3$, respectively.

\begin{theorem}
\label{thm:problem-1}
Let $G_1 = (V, E_1)$, $G_2 = (V, E_2)$, and $G_3 = (V, E_3)$ be three
simple graphs on the same vertex set $V$.
%
%
Let $\delta_1$ (resp.~$\Delta_2$, $\Delta_3$) denote the minimum
(resp.~maximum) degree of a vertex in $V(G_1)$ (resp.~$V(G_2)$,
$V(G_3)$). Let $\ell = o(\sqrt{n})$. Assume $\delta_1$, $\Delta_2$,
and $\Delta_3$ satisfy the following for $m=|V|$:
\[
m\left(\frac{(2\ell-1)^2-1}{(2\ell-1)^2}\right)+\frac{4\ell-3+q}{(2\ell-1)^2}\leq\delta_1\tag{1.1}\label{equation:condition-1}
\]
and
\[\Delta_2\Delta_3\leq\frac{q-1}{5\ell(2\ell-1)},\tag{1.2}\label{equation:condition-2}
\]
where $q\geq 1$ is an integer. Then there exists a Hamiltonian cycle $H$ in
$G_1$ such that $H^{\ell-1}$ (as a subgraph of ${G_1}$) does not form
an alternating-$(\ell,2,\ell,2)$-bag with $G_2$ and $G_3$.
\end{theorem}


This result follows a long line of work on finding Hamiltonian cycles and their powers in graphs with large minimum degree, which in particular includes a well-known theorem of
Dirac~\cite{dirac1952some} as well as a generalization of this
theorem~\cite{komlos1998posa}.

\begin{theorem}[Dirac~1952]
\label{thm:dirac}
Let $G$ be a simple graph on $m$ vertices. If the
minimum degree of a vertex in $G$ is at least $m/2$, then $G$ contains
a Hamiltonian cycle.
  \end{theorem}

In 1974, Seymour conjectured a more general result for Hamiltonian
cycles to the $p$th power \cite{seymour1974problem}, and in 1998,
Koml\'os, S\'ark\"ozy, and Szemer\'edi proved this
conjecture~\cite{komlos1998posa}.

\begin{theorem}[Koml\'os, S\'ark\"ozy, and Szemer\'edi 1998]\label{thm:dirac-type}
Let $G$ be a simple graph on $m$ vertices. If the minimum degree of a
vertex in $G$ is at least $m(p-1)/p$, then $G$ contains $H^{p-1}$
for a Hamiltonian cycle $H$.
\end{theorem}

Theorem~\ref{thm:problem-1} is a result of a similar form with
weakened conditions and a stronger conclusion for Hamiltonian cycles to
the $p$th power with a certain property, which we finally use to
obtain the following Baranyai theorem-related result that generalizes
Theorem \ref{thm:katona-katona-1}.

\begin{theorem}
\label{thm:baranyai-1}
Let $\alpha,\beta\in(0,1)$ such that $\alpha+\beta\geq 1$.  
Let $n$ be a positive integer. Let $k$ be $O(1)$ and $\ell$ be
$o(\sqrt{n})$. There exists a cyclic ordering of the $k$-element
subsets of $[n]$ with the following properties:

\begin{itemize}
\item[(1)]
any $\ell$ consecutive $k$-element subsets in the ordering are
disjoint, and in particular, form a $(k, \ell)$-parpartition of $[n]$,
which we refer to as a consecutive-$(k, \ell)$-parpartition; and
\item[(2)]
given any consecutive-$(k, \ell)$-parpartition $P_1$, no disjoint
consecutive-$(k, \ell)$-parpartition $P_2$ is too close to $P_1$
(that is, $P_1$ and $P_2$ are not $(\alpha,\beta)$-close).
\end{itemize}
\end{theorem}

\paragraph{Outline of the paper.}
We begin by introducing in Section~\ref{section:preliminaries} some
notations and definitions. In
Section~\ref{section:problem-3}, we prove
Theorem~\ref{thm:modified-katona}. In Section~\ref{section:problem-1}, we
prove Theorems~\ref{thm:problem-1} and~\ref{thm:baranyai-1}. Finally, we conclude with some open
questions.






  

\section{Preliminaries}
\label{section:preliminaries}
We briefly mention some notations on graphs in the following
presentation, which are largely standard.

Consider a graph $G=(V, E)$.
Given two vertices $u, v\in V$, we use $\edge{u}{v}$ for the edge in
$E$ connecting $u$ and $v$. Given vertices $v_0,v_1,\ldots,v_r$ for
$r\geq{0}$, we use $\path{v_0,v_1,\ldots,v_r}$ for the path consisting
of edges $\edge{v_i}{v_{i+1}}$ in $E$ for $0\leq{i} < r$, whose length
equals $r$.  We use $d_G(u, v)$ for the length of a shortest path in
$G$ from $u$ to $v$. We may write $d(u, v)$ for $d_G(u, v)$ if $G$ can
be readily inferred contextually.  Given a subset $S$ of $V$, we use
$G[S]$ for the subgraph of $G$ induced by the set of vertices $S\sub
V$. Formally, $V(G[S])=S$ and $E(G[S])$ consists of the edges in $E$
connecting vertices in $S$.

A Hamiltonian cycle $H$ in a graph $G$ is a subgraph of $G$ such that
$V(H)=V(G)$ and $E(H)$, a subset of $E(G)$, consists of the edges:
$$\edge{v_{0}}{v_{1}},\edge{v_{1}}{v_{2}},\ldots,\edge{v_{m-1}}{v_{m}}$$
for some enumeration of $v_0,v_1,\ldots,v_{m-1}$ of $V(G)$ (where
$n\geq 3$ and $v_{m}=v_0$).  In other words, $H$ is a circular path
that contains every vertex in $V(G)$ exactly once. We use
$H(v_0,v_1,\ldots,v_{m-1})$ for a Hamiltonian cycle consisting of
the circular path $\path{v_0,v_1,\ldots,v_{m-1},v_m}$ where $v_m=v_0$.

Given two graphs $J_1$ and $J_2$, we say that $J_1$ and $J_2$ are
disjoint if $V(J_1)$ and $V(J_2)$ are disjoint (as sets).
Given a graph $J$, we use $J^{p}$ for the $p$th power of $J$, where
$p$ ranges over natural numbers. Formally, $V(J^{p})=V(J)$ and
$E(J^{p})$ consists of edges $\edge{u}{v}$ for vertices $u,v\in{V(J)}$
satisfying $d_J(u,v)\leq p$. For instance, given a Hamiltonian cycle
$H$ (in some graph), any two vertices are connected in $H^p$, the
$p$th power $H$, if they are connected in $H$ by a path of length at
most $p$.  We use $H^{p}(\vv)$ for the $p$th power of a
Hamiltonian cycle $H(\vv)$, where
$\vv=(v_0,v_1,\ldots,v_{m-1})$ is an enumeration of vertices of the
underlying graph. Clearly, given a Hamiltonian cycle $H$ in $G$, its
$p$th power $H^p$ may or may not be a subgraph of $G$.






\def\hhh#1{#1}
\def\Hswap{{\it swap}}
\def\vvswap{{\it swap}}

\begin{definition}
Given a sequence $\vv$ of vertices $v_0, \ldots, v_{m-1}$ and an
integer $r\in[m]$, we use $\vv[i:r]$ for the subsequence
$v_i,v_{i+1},\ldots,v_{i+r-1}$ (of length $r$), where $0\leq i < m$
and arithmetic modulo $m$ is used.
\end{definition}

Assume $H(\vv)$ is a Hamiltonian cycle of $G$.  We use $H[i:r]$ for
the path in $H$ from $v_{i}$ to $v_{i+r-1}$.  We often use the phrase
{\em a copy of $H[r]$} to refer to $H[i:r]$ for some $i$ (which is
just a path of length $r-1$ in $H$).  Similarly, we use $H^{p}[i:r]$ for
$H^{p}[\vv[i:r]]$, that is, the subgraph obtained from restricting
$H^{p}$, the $p$th power of $H$, onto the vertex set of $\vv[i:r]$,
and refer to $H^{p}[i:r]$ as a copy of $H^{p}[r]$ for each $0\leq{i}<m$.




Given a graph $J$, we use $K_\ell(J)$ to refer to a subgraph of $J$
that is isomorphic to $K_\ell$. And $K_{\ell,\ell}(J)$ is similarly
interpreted.

Given $H=H(\vv)$ for $\vv=(v_0,v_1,\ldots,v_{m-1})$, each $H^p[i:p+1]$
for $0\leq{i}<m$ is a copy of $K_{p+1}$, since by definition of $H^p$,
any two vertices no further than $p$ edges apart along the cycle are
adjacent; and as these are the only adjacent pairs of vertices, $H^p$
contains exactly $m$ copies of $K_{p+1}$.


Let $G_1 = (V,E_1)$, $G_2 = (V,E_2)$ and $G_3=(V,E_3)$ be three simple
graphs such that $E_1\cap{E_2}=E_1\cap{E_3}=\emptyset$.  Let $J_1$ be
a (not necessarily proper) subset of $G_1$. We define the following
relationship between copies of $K_\ell$ in $J_1$ connected by the
edges in $G_2$ and $G_3$.

\begin{definition}
We write that $J_1$ forms an alternating-$(\ell, 2,
\ell, 2)$-bag with $G_2$ and $G_3$ if there exist two copies $A_1$ and $A_2$ of $K_\ell$ in
$J_1$ and two vertex-disjoint edges $e_2\in E(G_2)$ and $e_3\in
E(G_3)$ such that $e_2$ and $e_3$ each connect a vertex in $A_1$ to a
vertex in $A_2$.
\end{definition}
We may also write that $J_1$ forms an alternating-$(\ell, 2, \ell,
2)$-bag with $G_2$ and $G_3$ via $A_1$ and $A_2$, or simply that $A_1$
and $A_2$ form an alternating-$(\ell, 2, \ell, 2)$-bag if it is clear
from the context what $J_1$, $G_2$, and $G_3$ are.  Note that this
definition of an alternating-$(\ell, 2, \ell, 2)$-bag is referred to
as the three-graph version (in contrast to the two-graph version given
earlier).

\section{Existence of linear-sized ``not-too-close'' $(k, \ell)$-parpartitions of $[n]$}
\label{section:problem-3}

In this section, we prove Theorem~\ref{thm:modified-katona}, which
we deduce from the following theorem on graphs.
\begin{theorem}\label{thm:modified-graphs}
Let $G_1 = (V, E_1)$, $G_2 = (V, E_2)$, and $G_3 = (V, E_3)$ be three
simple graphs on the same vertex set $V$, where $|V| = m$. Let
$\delta_1$ (resp.~$\Delta_2$, $\Delta_3$) denote the minimum
(resp.~maximum) degree of a vertex in $V(G_1)$ (resp.~$V(G_2)$,
$V(G_3)$). Let $\ell\geq 2$ be an integer. Assume $\delta_1$,
$\Delta_2$, and $\Delta_3$ satisfy
\[m\left(\frac{3\ell-1}{3\ell}\right)\leq\delta_1\tag{2.1}\label{equation:condition-3}
\]
and
\[\Delta_2\Delta_3\leq\frac{m-3}{15\ell^2},\tag{2.2}\label{equation:condition-4}.
\]
Then there are $\lfloor m/\ell\rfloor$ pairwise disjoint copies of $K_\ell$
in $G_1$ such that no two form an alternating-$(\ell, 2, \ell,
2)$-bag with $G_2$ and $G_3$.
\end{theorem}

First, we show that if Condition~\ref{equation:condition-3} holds, then $G_1$ contains $\lfloor
m/\ell\rfloor$ pairwise disjoint copies of $K_\ell$, which we refer to as an \emph{almost-$\ell$-decomposition} of $G_1$.

\begin{lemma}\label{lemma:modified-complete}

Let $G = (V, E)$ be a simple graph on $m$ vertices. Let $\delta$
be the minimum degree of a vertex in $G$, and assume that 
\[m\left(\frac{3\ell-1}{3\ell}\right)\leq\delta.
\]
Then $G$ has an almost-$\ell$-decomposition.

\end{lemma}

\begin{proof}

Let $G_0 = G$ as given. If $G = K_m$, then $G$ clearly has an
almost-$\ell$-decomposition. Suppose $\delta$ satisfies the given
condition and $G$ does not have an almost-$\ell$-decomposition, and
let us add edges to $G$ until it does. Let us call this
almost-$\ell$-decomposition $L$. Let $e = (a_1, a_2)$ be the last edge
added, and let $G' = (V, E' = E\backslash\{e\})$ be $G$ just prior to
adding $e$, so that $G'$ contains $L$ with one edge missing from one
copy of $K_\ell$, which we denote by $A_1$. We show that $G'$ has an
almost-$\ell$-decomposition, which would be a contradiction, implying
that $G_0$ does have an almost-$\ell$-decomposition. We use the
following claim, which identifies a set of ``swapping candidate" vertices $b_1\notin V(A_1)$. 

\begin{claim}\label{claim:good-vertices}

There are at least $m/3$ vertices $b_1\notin V(A_1)$ such that (1) $(a, b_1)\in E(G')$ for all $a\in A_1$
and (2) $(a_1, b)\in E'$ for all $b\in B_1$, where $B_1$ is the copy
of $K_\ell$ in $L$ whose vertex set contains $b_1$.
  
\end{claim}

\begin{proof}[Proof of Claim~\ref{claim:good-vertices}]

Define a \emph{swapping candidate vertex} to be one that satisfies the two properties in the statement of the claim and a \emph{swapping noncandidate vertex} to be one that does not satisfy both properties. Furthermore, define a \emph{swapping noncandidate of type 1 (resp. 2)} to be a vertex that does not satisfy property 1 (resp. 2). We show as follows that the number of swapping noncandidates is strictly less than the number of swapping candidates minus $m/3$. Let $\overline{G'} = (V, \overline{E'})$.

If $b_1$ is a swapping noncandidate of type 1, then there is an edge of the form $(a, b_1)$, where $a\in A_1$, in $\overline{E'}$, which we think of as witnessing this. Each edge of the form $(a, v)$, where $v\in V$, can witness at most one swapping noncandidate of type 1; thus $\#\{(a, v)\in\overline{E'}\mid a\in A_1\text{ and }v\in
V\}\leq\ell(m-\delta)$ upper bounds the number of swapping noncandidates of type 1. 

If $b_1$ is a swapping noncandidate of type 2, then $b_1$ resides in a copy of $B_1$ of $K_\ell$ in $L$ for which there exists a vertex $b_2\in B_1$ such that $(a_1, b_2)\notin E$; we refer to $B_1$ as a "bad" copy of $K_\ell$. If $B_1$ is bad, then there is an edge of the form $(a_1, b)\in\overline{E'}$, where $b\in B_1$, that we think of as witnessing this. Each edge of the form $(a_1, v)$, where $v\in V$, witnesses at most one bad copy of $K_\ell$ because of disjointness. Thus, 
\[
\#\{(a_1, v)\in\overline{E'}\mid v\in V\backslash V(A_1)\}\leq m-\delta-(|A_1|-1)\leq m-\delta-1
\]
is an upper bound on the number of bad copies of $K_\ell$ in $L$.

Now, in order to show that there are at least $m/3$ swapping candidates, it suffices to show that 
\[\ell(m-\delta)+\ell(m-\delta-1)+\frac{m}{3} < m.
\]
For this, we need only to show that $m(3\ell-1)/3\ell-1/2 < \delta$,
which follows from the assumption that 
$m(3\ell-1)/3\ell\leq\delta$.
\end{proof}
We can obtain an almost $\ell$-decomposition of $G'$ by swapping $a_1$ with any swapping candidate vertex $b_1$.
\end{proof}

Next, we show that $G_1$ has an almost-$\ell$-decomposition $L$ such that no two copies of of $K_\ell$ in $L$ form an alternating-$(\ell, 2, \ell, 2)$-bag with $G_2$ and $G_3$. 

\begin{lemma}\label{lemma:modified-no-bag}
  
Let $G_1 = (V, E_1)$, $G_2 = (V, E_2)$, and $G_3 = (V, E_3)$ be three
graphs that satisfy the conditions of
Theorem~\ref{thm:modified-graphs}. Let $A_1$ and $A_2$ be two copies
of $K_\ell$ in an almost-$\ell$-decomposition $L$ of $G_1$. Let $(a_1, a_2)$
be an edge, where $a_1\in A_1$ and $a_2\in A_2$, such that $(a_1,
a_2)\notin E_1\cup E_3$.

Assume that no two copies of $K_\ell$ in $L$ form an
alternating-$(\ell, 2, \ell, 2)$-bag. Then there exists a vertex $b_1$ not in
$L$ such that swapping $a_1$ and $b_1$ in $G_1$ results in a new
almost-$\ell$-decomposition $L'$ of $G_1$ such that no two copies of $K_\ell$
in $L'$ form an alternating $(\ell, 2, \ell, 2)$-bag, even if the edge
$(a_1, a_2)$ is added to $E_3$.
  
\end{lemma}

\begin{proof}

By Claim~\ref{claim:good-vertices} (stated in the proof of
Lemma~\ref{lemma:modified-complete}), there are at least $m/3$ swapping candidate
vertices $b_1$ in $V$ such that $(a, b_1)\in E_1$ for all $a\in A_1$
and $(a_1, b)\in E_1$ for all $b\in B_1$, where $B_1$ is the copy of
$K_\ell$ in $L$ whose vertex set contains $b_1$. If we swap $a_1$ with any
swapping candidate vertex in the almost-$\ell$-decomposition $L$, we obtain a new
almost-$\ell$-decomposition $L'$ of $G_1$. So it suffices to show that the
number of swapping candidates for which two copies of $K_\ell$ in $L'$
form an alternating-$(\ell, 2, \ell, 2)$-bag, which we refer to as \emph{bad
swapping candidates}, is strictly less than the total number of swapping candidates. Let
$b_1$ be a swapping candidate in a copy $B_1$ of $K_\ell$ in $L$, and we consider as follows three
cases where $b_1$ is a bad candidate. 

First, suppose that we have
swapped $a_1$ and $b_1$ in $L$, obtaining a new almost-$\ell$-decomposition $L'$, and have not yet added the edge $(a_1, a_2)$ to $E_3$. Recall that we assume that the edges in $G_1$, $G_2$, and $G_3$ are colored red, green, and blue, respectively.

\begin{itemize}

\item[(1)] Case 1: There is a clique $D\neq A_1$ in $L'$ such that $A_1$ forms an alternating-$(\ell, 2, \ell, 2)$-bag with $D$. In
  this case, there exists either a blue-red-green path or a
  green-red-blue path $(a_1, d, d', b)$, where $d, d'\in V(D)$ and
  $b\in B_1$, which we think of as a witness to $b_1$ being a bad
  swapping candidate. Each blue-red-green and green-red-blue path of the
  form $(a_1, d, d', v)$, where $d, d'\in V(D)$ and $v\in V$, can act
  as a witness for at most $\ell$ bad swapping candidates. Starting from $a_1$,
  for a blue-red-green path (resp.~green-red-blue path), there are at
  most $\Delta_2$ (resp.~$\Delta_3$) choices for $d$, $\ell$ choices
  for $d'$, and $\Delta_3$ (resp.~$\Delta_3$) choices for $v$; hence,
  there are at most $2\ell\Delta_2\Delta_3$ witness paths, and
  $2\ell^2\Delta_2\Delta_3$ is an upper bound on the number of bad swapping candidates in this case.

\item[(2)] Case 2: There is a clique $D\neq B_1$ in $L'$ such that $B_1$ forms an alternating-$(\ell, 2, \ell, 2)$-bag with $D$. In
  this case, there exists either a blue-red-green path or a
  green-red-blue path $(a, d, d', b_1)$, where $a\in A_1$ and $d, d'\in
  V(D)$, which we think of as a witness to $b_1$ being a bad swapping candidate. Each blue-red-green path and green-red-blue path of the
  form $(a, d, d', v)$, where $a\in A_1$, $d, d'\in V(D)$, and $v\in
  V$, can act as a witness for at most one bad swapping candidate. For a
  blue-red-green path (resp.~green-red-blue path), there are at most
  $\ell$ choices for $a$, $\Delta_2$ (resp.~$\Delta_3$) choices for
  $d$, $\ell$ choices for $d'$, and $\Delta_3$ (resp.~$\Delta_2$)
  choices for $v$; hence, there are at most $2\ell^2\Delta_2\Delta_3$
  witness paths, and $2\ell^2\Delta_2\Delta_3$ is an upper bound on
  the number of bad swapping candidates in this case.

\end{itemize}

Now suppose that we add the edge $(a_1, a_2)$ to $E_3$. We consider one
more case.

\begin{itemize}

\item[(3)] Case 3: $B_1$ forms an alternating-$(\ell, 2, \ell,
  2)$ bag with $A_2$. In this case, there exists a blue-red-green path $(a_1,
  a_2, a', b)$, where $a'\in A_2$ and $b\in B_1$, which we think of as
  a witness to $b_1$ being a bad swapping candidate. Each blue-red-green
  path of the form $(a_1, a_2, a', v)$, where $a'\in A_2$ and $v\in
  V$, can act as a witness for at most $\ell$ bad swapping candidates. There
  are at most $\ell$ choices for $a'$ and $\Delta_2$ choices for $v$;
  hence, there are at most $\ell\Delta_2$ witness paths, and
  $\ell^2\Delta_2$ is an upper bound on the number of bad swapping candidates in this case.
  
\end{itemize}

Summing over all of the cases, we obtain that the total number of bad
swapping candidates is at most $2\ell^2\Delta_2\Delta_3 +
2\ell^2\Delta_2\Delta_3 + \ell^2\Delta_2\leq
5\ell^2\Delta_2\Delta_3$. Thus, it suffices to show the following:
\[5\ell^2\Delta_2\Delta_3+1\leq\frac{m}{3}, 
\]
which is the same as
\[\Delta_2\Delta_3\leq\frac{m-3}{15\ell^2},
\]
which is assumed.
\end{proof}

Finally, we use Theorem~\ref{thm:modified-graphs} to prove Theorem~\ref{thm:modified-katona}.

\begin{proof}[Proof of Theorem~\ref{thm:modified-katona}]
We use Theorem~\ref{thm:modified-graphs}. Let $V = {[n]\choose
  k}$. Let two vertices be adjacent in $G_1$ if
their corresponding $k$-element subsets are disjoint. Let two vertices be adjacent in $G_2$ (resp.~$G_3$) if
their corresponding $k$-element subsets intersect in more than $\alpha k$ (resp.~$\beta k$) elements. We have that a $(k,
\ell)$-parpartition corresponds to a copy of $K_\ell$ in $G_1$ and
that two parpartitions are $(\alpha, \beta)$-close if and only if the corresponding cliques in $G_1$
form an alternating-$(\ell, 2, \ell, 2)$-bag with $G_1$ and $G_2$. Hence, in order to prove the theorem, it suffices
to show that $G_1$, $G_2$, and $G_3$ satisfy the degree conditions of
Theorem~\ref{thm:modified-graphs}.

We have that $G_1$ is regular and that the degree of each vertex is
${n-k\choose k}$. Condition~\ref{equation:condition-3}, then, is the following:
\[{n\choose k}\left(\frac{3\ell-1}{3\ell}\right)\leq {n-k\choose k}.
\]
This inequality holds when $k$ is a constant and $n$ is large if
$k^2/n\leq 1/3\ell$, which is the same as $k^2\ell\leq n/3$, which is
assumed.

We have that $G_2$ and $G_3$ are also regular. The degree of each
vertex in $G_2$ is
\[\sum_{1\leq i < (1-\alpha)k} {k\choose i}{n-k\choose i},
\]
and the degree of each vertex in $G_3$ is
\[\sum_{1\leq i < (1-\beta)k} {k\choose i}{n-k\choose i}.
\]
Then, Condition~\ref{equation:condition-4} states that, for $m={n\choose k}$,
\[\left(\sum_{1\leq i < (1-\alpha)k} {k\choose i}{n-k\choose i}\right)\left(\sum_{1\leq i < (1-\beta)k} {k\choose i}{n-k\choose i}\right)\leq\frac{m-3}{15\ell^2},
\]
which holds asymptotically as the left hand side is
$O(n^{(2-(\alpha+\beta))k-1})$, which is $O(n^{k-3})$ when
$\alpha+\beta\geq (k+2)/k$ and the right hand side is
$\omega(n^{k-2})$.

Hence, as $G_1$, $G_2$, and $G_3$ satisfy both Conditions \ref{equation:condition-3} and \ref{equation:condition-4}, we conclude that the number of $(k, \ell)$-parpartitions of $n$
such that no $k$-element subset appears in two parpartitions and no
two parpartitions are too close is $\lfloor{m/\ell}\rfloor$ for
$m={n\choose k}$.
\end{proof}  


\section{Existence of $(\ell, 2, \ell, 2)$-bag-avoiding powers of Hamiltonian cycles}
\label{section:problem-1}
In this section, we prove Theorem~\ref{thm:baranyai-1}. To do this, we first
prove a more general graph result in Theorem~\ref{thm:problem-1}, which we restate here.

\vspace{6pt}
\noindent{\bf Theorem~\ref{thm:problem-1}.}~
{\em
Let $G_1 = (V, E_1)$, $G_2 = (V, E_2)$, and $G_3 = (V, E_3)$ be three
simple graphs on the same vertex set $V$.
%
%
Let $\delta_1$ (resp.~$\Delta_2$, $\Delta_3$) denote the minimum
(resp.~maximum) degree of a vertex in $V(G_1)$ (resp.~$V(G_2)$,
$V(G_3)$). Let $\ell = o(\sqrt{n})$. Assume $\delta_1$, $\Delta_2$,
and $\Delta_3$ satisfy the following for $m=|V|$:
\[
m\left(\frac{(2\ell-1)^2-1}{(2\ell-1)^2}\right)+\frac{4\ell-3+q}{(2\ell-1)^2}\leq\delta_1\tag{1.1}\label{equation:condition-1}
\]
and
\[\Delta_2\Delta_3\leq\frac{q-1}{5\ell(2\ell-1)},\tag{1.2}\label{equation:condition-2}
\]
where $q\geq 1$ is an integer. Then there exists a Hamiltonian cycle $H$ in
$G_1$ such that $H^{\ell-1}$ (as a subgraph of ${G_1}$) does not form
an alternating-$(\ell,2,\ell,2)$-bag with $G_2$ and $G_3$.
}

In this first lemma, we show that if
Condition~\ref{equation:condition-1} holds, then $G_1$ contains a
$H^{\ell-1}$ for some Hamiltonian cycle $H$.

\begin{lemma}\label{lemma:ham-cycle}
Let $G = (V, E)$ be a simple graph on $m$ vertices with minimum degree $\delta = \delta(G)$. Suppose $\delta$ satisfies
\[m\left(\frac{(2\ell-1)^2-1}{(2\ell-1)^2}\right)+\frac{4\ell-3+q}{(2\ell-1)^2}\leq\delta.
\]
Then $G$ contains $H^{\ell-1}$ for some Hamiltonian cycle $H$.
\end{lemma}  

\begin{proof}
%
%
Let $G_0 = G$ as given. If $G = K_m$, then the lemma clearly holds.
Suppose that $G$ does not contain $H^{\ell-1}$ for any Hamiltonian
cycle $H$.  Let us keep adding edges to $E$ until $G$ contains
$H_0^{\ell-1}$ for some Hamiltonian cycle $H_0$. Without loss of generality, we assume
$H_0=H(v_0, \ldots, v_{m-1})$ is contained in $G_0$. Let $e_1$ be the
last edge added to $E$ and $G' = (V, E' = E(G)\backslash\{e_1\})$ be
the graph just prior to adding $e_1$. Hence, $G'$ contains
$H^{\ell-1}_0$ with one edge $e_1$ missing, which we denote by
$\hat{H}^{\ell-1}_0$
(that is, $\hat{H}^{\ell-1}_0 = {H}^{\ell-1}_0\backslash\{e_1\}$).

We show that $G'$ contains $H^{\ell-1}$ for some $H$. As this would be
a contradiction, it implies that $G_0$ does contain $H^{\ell-1}$ for
some $H$. We make use of the following claim.

\begin{claim}\label{claim:span-bipartite}
Given $C_0=\hhh{H}_0[i_0:2\ell-1]$ for some $i_0$, there are at least
$q$ copies of $\hhh{H}_0[2\ell-1]$ disjoint from $C_0$ that form
$K_{2\ell-1, 2\ell-1}(G')$ with $C_0$, where a copy of
$\hhh{H}_0[2\ell-1]$ refers to $\hhh{H}_0[i:2\ell-1]$ for some $i$.
\end{claim}

\begin{proof}[Proof of Claim~\ref{claim:span-bipartite}]
Let $\overline{G'} = (V, \overline{E'})$ be the complement of $G'$.  Let
$W_0$ be the vertex set of $C_0$, and let $W$ be the vertex set of a copy of $\hhh{H}_0[2\ell-1]$ disjoint from
$C_0$. We have that $W_0$ and $W$ do not span $K_{2\ell-1,
  2\ell-1}(G')$ if and only if there exists an edge
$\edge{u}{v}\in\overline{E'}$ such that $u\in W_0$ and $v\in W$, which
we think of as a \emph{witness} to $W$ being ``bad''; $W$ is bad if and only if there
exists a witness edge for $W$. A particular edge $\edge{u}{v}$, where $u\in
W_0$ and $v\in V$, can act as a witness for at most $2\ell-1$ bad
copies of $\hhh{H}_0[2\ell-1]$, since $v$ is in at most $2\ell-1$
copies of $\hhh{H}_0[2\ell-1]$. Hence, $(2\ell-1)\#\{(u, v)\mid u\in
W_0\text{ and }v\in V\}$ is an upper bound on the
number of bad copies of $\hhh{H}_0[2\ell-1]$, which do not
span $K_{2\ell-1, 2\ell-1}(G')$ with $W_0$.

There are $m$ distinct copies of $\hhh{H}_0[2\ell-1]$. Out of
them, $m-(4\ell-3)$ copies are disjoint from $C_0$. To show the
claim, it suffices to show that the sum of the number of bad copies of
$\hhh{H}_0[2\ell-1]$ and $q$, which we think of as the number
of good copies of $\hhh{H}_0[2\ell-1]$, is at most the total
number of copies of $\hhh{H}_0[2\ell-1]$ disjoint from $C_0$. In other words, it suffices to show that
\[(2\ell-1)\#\{(u, v)\mid u\in W_0\text{ and }v\in V\}+q\leq m-(4\ell-3).
\]
We have $\#\{(u, v)\mid u\in W_0\text{ and }v\in
V\}\leq (2\ell-1)(m-\delta)$, since $|W_0| = 2\ell-1$
and the maximum degree of a vertex in $\overline{G'}$ is $m-\delta$
(as the minimum degree of a vertex in $G'$ is $\delta$). Hence, it
suffices to show
\[\delta\geq m\left(\frac{(2\ell-1)^2-1}{(2\ell-1)^2}\right)+\frac{4\ell-3+q}{(2\ell-1)^2},
\]
which is assumed. 
\end{proof}

\begin{claim}\label{claim:cor-bipartite}
Let $r\in [2\ell-2]$. Given a copy $C_1$ of $\hhh{H}_0[r]$, there is a
copy $C_2$ of $\hhh{H}_0[r]$ disjoint from $C_1$ such that $V(C_1)$
and $V(C_2)$ span $K_{r,r}(G')$.
\end{claim}

\begin{proof}[Proof of Claim~\ref{claim:cor-bipartite}]
Let $C'_1$ be a copy of $\hhh{H}_0[2\ell-1]$ containing $C_1$ as a
subgraph. By Claim~\ref{claim:span-bipartite}, there exists a copy
$C'_2$ of $\hhh{H}_0[2\ell-1]$ disjoint from $C'_1$ with vertex set
$W_2'$ such that $W_1'$ and $W_2'$ span $K_{2\ell-1, 2\ell-1}$. Let
$C_2$ be a copy of $\hhh{H}_0[r]$ that is a subgraph of $C'_2$, and we
have that $C_1$ and $C_2$ are disjoint and they span $K_{r,r}(G')$.
\end{proof}  

\begin
{figure}[t]
\centering
\begin
{tikzpicture}
[scale=1.00]
[line cap=rect,line width=4.5pt]
\filldraw%
[fill=white]%
(0,0)circle[radius=2.5cm];
\foreach\angle[count=\xi]in{60,30,...,-270}
{
\filldraw[fill=black]%
  (\angle:2.50cm)circle[radius=1.75pt];
\node[font=\large]at(\angle:2.85cm){\mbox{$v_{\xi}$}};
\draw(060:2.5cm)--(000:2.5cm); 
\draw(000:2.5cm)--(300:2.5cm); 
\draw(300:2.5cm)--(240:2.5cm); 
\draw(240:2.5cm)--(180:2.5cm); 
\draw(180:2.5cm)--(120:2.5cm); 
\draw(120:2.5cm)--(060:2.5cm); 
\draw(030:2.5cm)--(330:2.5cm); 
\draw(330:2.5cm)--(270:2.5cm); 
\draw(270:2.5cm)--(210:2.5cm); 
\draw(210:2.5cm)--(150:2.5cm); 
\draw(150:2.5cm)--(090:2.5cm); 
\draw[red][dashed](090:2.5cm)--(030:2.5cm); 
\draw[blue][dotted](060:2.5cm)--(210:2.5cm); 
\draw[blue][dotted](090:2.5cm)--(240:2.5cm); 
\draw[blue](030:2.50)arc(030:060:2.50); 
\draw[blue](000:2.50)arc(000:030:2.50); 
\draw[blue](330:2.50)arc(330:360:2.50); 
\draw[blue](300:2.50)arc(300:330:2.50); 
\draw[blue](270:2.50)arc(270:300:2.50); 
\draw[blue](240:2.50)arc(240:270:2.50); 
\draw[blue](090:2.50)arc(090:120:2.50); 
\draw[blue](120:2.50)arc(120:150:2.50); 
\draw[blue](150:2.50)arc(150:180:2.50); 
\draw[blue](180:2.50)arc(180:210:2.50); 
}
\end{tikzpicture}
\caption{%
A drawing of some $H^2$ with one
edge, colored red (dashed), missing. The rest of the edges in $H^2$
are solid lines, colored either blue or black.
The cycle in blue (consisting of both dotted and solid lines)
is a new Hamiltonian cycle $H'$ such that $H'^2$ is (fully) contained
in the graph $G$; to show this, we consider a copy of $H[4]$ surrounding
the missing edge (here consisting of the vertices $v_{11}, v_{12}, v_1, v_2$)
and another disjoint copy of $H[4]$ (here consisting of the vertices
$v_6, v_7, v_8, v_9$), which together form $K_{4, 4}$. (The existence of
the other disjoint copy of $H[4]$ is given by Claim~\ref{claim:cor-bipartite}.)
For instance, the two dotted blue lines refer to two edges in this $K_{4,4}$.
Essentially, because of this $K_{4, 4}$, in the $H'^2$ induced by the new ordering
$H'$, the vertices that have new neighbors (that is, neighbors
different from theirs in $H^2$) are indeed connected by an edge to
each of these new neighbors.}
\label{figure:applying_corollary_bipartite}
\end{figure}
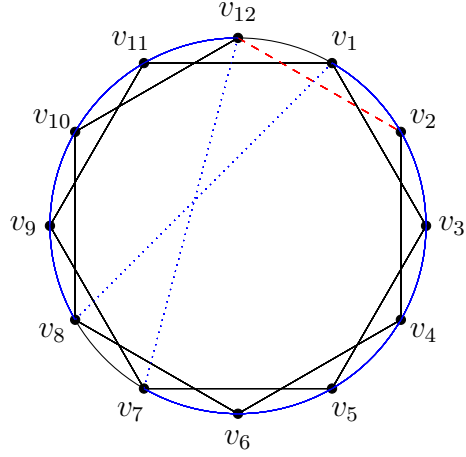

We now use Claim~\ref{claim:cor-bipartite} to show that $G'$ contains
$H^{\ell-1}$ as a subgraph for some Hamiltonian cycle $H$.  In the
following, arithmetic operations on vertex subscripts are modulo $m$.

We have that $G'$ contains $H^{\ell-1}_0$ with the edge $e_1$ missing.
Write $e_1=\edge{v_i}{v_{i'}}$ with $i < i'$ and set $d = i'-i$; we have that $1\leq d\leq \ell-1$.
Let $S$ be the sequence
$$\{v_{i-\ell+2},\ldots,v_i,\ldots,v_{i'},\ldots,v_{i+\ell-1}\},$$
which is a segment of $v_0,v_1,\ldots,v_{m-1}$ of length
$2\ell-2$.  Consider the subgraph $\hhh{H}_0[S]$, which is a copy of
$\hhh{H}_0[2\ell-2]$. By Claim~\ref{claim:cor-bipartite},
there exists another copy $C^*$ of $\hhh{H}_0[2\ell-2]$ disjoint from
$\hhh{H}_0[S]$ such that $W^*$ and $S$ span $K_{2\ell-2, 2\ell-2}(G')$
for $W^*=V(C^*)$.  We have $W^* =
\{v_{j},v_{j+1},\ldots,v_{j'-1},v_{j'}\}$ for some $j < j'$ (which
implies $j'=j+2\ell-3$). Essentially, in order to obtain
$H_1^{\ell-1}$ in $G'$ for some Hamiltonian cycle $H_1$, we rearrange
the vertices in $W^*$ and $S$ in $\hhh{H}_0$, and the corollary
ensures that every pair of vertices $u, v\in V(G')$ is connected by an
edge as long as $d(u, v)\leq p$ in ${H_1}$. We define $H_1$
as the Hamiltonian cycle
\[
H_1 = (v_{i+1}, \ldots, v_j, \ldots, v_{j+\ell-2}, v_{i}, v_{i-1}, \ldots, v_{j'}, v_{j'-1}, \ldots, v_{j'-\ell+2}).
\]
In other words, $H_1$ takes the same path as $H_0$ from $v_{i+1}$ to $v_{j+\ell-2}$; it then
jumps from $v_{j+\ell-2}$ to $v_{i}$ and takes a path in $H_0$ from $v_{i}$ to $v_{j'-\ell+2}$;
it finally closes the cycle with a jump from $v_{j'-\ell+2}$ to $v_{i+1}$. For instance, in Figure~\ref{figure:applying_corollary_bipartite}, $H_1$ is the cycle in blue, where $v_{i+1} = v_1$, $v_{j+\ell-2} = v_7$, $v_i = v_{12}$, and $v_{j'-\ell+2} = v_8$.  
To show that $H_1^{\ell-1}$ is indeed a Hamiltonian cycle to the
$(\ell-1)$th power in $G'$, it suffices to show that each pair of vertices
that are neighbors in $H_1^{\ell-1}$ but not in $H_0^{\ell-1}$ is
connected by an edge in $E(G')$. The following vertices have the following new neighbors in
$H_1^{\ell-1}$:

\begin{itemize}
\item $v_j, \ldots, v_{j+\ell-2}$, with vertices in $\{v_i, v_{i-1}, \ldots, v_{i-\ell+2}\}$;
\item $v_{i+1}, \ldots, v_{i+\ell-1}$, with vertices in $\{v_{j'}, v_{j'-1},\ldots, v_{j'-\ell+2}\}$;
\item $v_i, v_{i-1}, \ldots, v_{i-\ell+2}$, with vertices in $\{v_j, \ldots, v_{j+\ell-2}\}$; and
\item $v_{j'}, v_{j'-1}, \ldots, v_{j'-\ell+2}$, with vertices in $\{v_{i+1}, \ldots, v_{i+\ell-1}\}$.   
\end{itemize}  

The vertices with new neighbors in $H_1^{\ell-1}$ are exactly those in
$S\cup W^*$; the vertices in $S$ have new neighbors in $W^*$ and the
vertices in $W^*$ have new neighbors in $S$. Since $S$ and $W^*$ span
$K_{2\ell-2, 2\ell-2}(G')$ (by
Claim~\ref{claim:cor-bipartite}), each pair of new neighbors
is connected by an edge in $E(G')$. Therefore, $G'$ contains
$H_1^{\ell-1}$, a Hamiltonian cycle raised to power $(\ell-1)$, contradicting our initial assumption.
\end{proof}

We next show that $G_1$ contains an $H^{\ell-1}$ for some Hamiltonian
cycle $H$ that does not form an alternating-$(\ell, 2, \ell, 2)$-bag
with $G_2$ and $G_3$.

\begin{lemma}\label{lemma:no-bag}
Let $G_1 = (V, E_1)$, $G_2 = (V, E_2)$, and $G_3 = (V, E_3)$ be three
graphs whose minimum/maximum degrees satisfy the inequalities (1.1) and (1.2) in Theorem~\ref{thm:problem-1}.

Let $H_1$ be a Hamiltonian cycle such that $H_1^{\ell-1}$ is in $G_1$. Let $A_1$ and
$A_2$ be two disjoint copies of $K_\ell$ in $H_1^{\ell-1}$. Let
$\edge{a_1}{a_2}$ be an edge, where $a_1\in A_1$ and $a_2\in A_2$, such
that $(a_1, a_2)\notin{E_3}$.

Assume that $H_1^{\ell-1}$ does not form an alternating-$(\ell,
2,\ell, 2)$-bag with $G_2$ and $G_3$. Then there exists a vertex
$a_3\in V$ such that swapping $a_1$ and $a_3$ in $H_1$ creates a new
Hamiltonian cycle $H_2$, where $H_2^{\ell-1}$ is in $G_1$ and
$H^{\ell-1}_2$ does not form an alternating-$(\ell, 2, \ell, 2)$-bag
with $G_2$ and $G_3\cup\{\edge{a_1}{a_2}\}$.
\end{lemma}

\begin{proof}
Write $H_1=(v_0,v_1,\ldots, v_{m-1})$. Write $a_1 = v_{i_1}$ and $a_2 = v_{i_2}$ for $0\leq i_1, i_2\leq
m-1$. For each $j = 1, 2$, let $C_j$ be the subgraph of $H_1^{\ell-1}$
induced by $\{v_{i_j-\ell+1}, \ldots, v_{i_j}, \ldots,
v_{i_j+\ell-1}\}$. For $j = 1, 2$, let $D_{j, 1}, \ldots, D_{j, \ell}$ be the $\ell$
copies of $K_\ell$ in $H_1^{\ell-1}$ with vertex sets containing
$v_{i_j}$. Note that each $C_j$ is a copy of $H_1^{\ell-1}[2\ell-1]$
and $C_j = \bigcup_{l\in [\ell]} D_{j, l}$.

By Claim~\ref{claim:span-bipartite} (stated in the proof of
Lemma~\ref{lemma:ham-cycle}), there are at least $q$ copies of
$H_1[2\ell-1]$ disjoint from $C_1$, whose vertex sets we denote by
$U_1, \ldots, U_q$, such that for each $r\in[q]$, $V(C_1)$ and $U_r$
span $K_{2\ell-1, 2\ell-1}(G_1)$. If $U_r = \{v_{j_r}, \ldots,
v_{j_r+2\ell-2}\}$, then let $c_r = v_{j_r+\ell-1}$, i.e., $c_r$ is
the ``center'' vertex in $V_r$. By swapping $a_1$ with any $c_r$ in
$H_1$, we obtain a new Hamiltonian cycle $\hat{H}_1$ such that
$\hat{H}_1^{\ell-1}$ is in $G_1$, since $a_1$ ($c_r$) is adjacent to
every vertex that $c_r$ ($a_1$) is adjacent to in $H_1^{\ell-1}$. To
show the lemma, it suffices to show that the number of center vertices
$c_r$ for which swapping $a_1$ and $c_r$ in $H_1^{\ell-1}$ results in
$\hat{H}_1^{\ell-1}$ forming an alternating-$(\ell, 2, \ell, 2)$-bag
with $G_2$ and $G_3$, which we refer to as ``bad'' $c_r$'s, is less
than $m$, which lower bounds the total number of $c_r$'s.

For $r\in [q]$, suppose $v_{i_r} = c_r$ for some $0\leq i_r\leq
m-1$. Let $T_r$ be the subgraph of $H_1^{\ell-1}$ induced by
$\{v_{i_r-\ell+1}, \ldots, v_{i_r}, \ldots, v_{i_r+\ell-1}\}$. Let
$R_{r,1}, \ldots, R_{r,\ell}$ be the $\ell$ copies of $K_\ell$ in
$H_1^{\ell-1}$ whose vertex sets contain $v_{i_r}$. Note that $T_r$ is
a copy of $H_1^{\ell-1}[2\ell-1]$ and $T_r = \bigcup_{l\in [\ell]}
R_{r, l}$. Let $c = c_{r^*}$, which is $v_{i_{r^*}}$ for some $r^*\in [m]$. We
consider as follows three cases where $c$ is bad. First, suppose we
swap $a_1$ with $c$ in $H_1^{\ell-1}$ and do not yet add $(a_1, a_2)$
to $E_3$.

\begin{itemize}

\item[(1)]

Case 1: $\hat{H}_1^{\ell-1}$ forms with $G_2$ and $G_3$ an
alternating-$(\ell, 2, \ell, 2)$-bag via $R = R_{{r^*},{l^*}}$ for
some $l^*\in [\ell]$ and some copy $A^*$ of $K_\ell(G_1)$ disjoint
from $R$.

If $R$ forms an alternating-$(\ell, 2, \ell, 2)$-bag with $A^*$, then
$T_{r^*}$ forms an alternating-$(2\ell-1, 2, \ell, 2)$-bag with $A^*$,
which implies the existence of either a blue-red-green path or a
green-red-blue path $(a_1, h, h', u)$, where $h, h'\in V(A^*)$ and
$u\in V(T_{r^*})$ that we think of as a witness to $c$ being
bad. Each blue-red-green path and green-red-blue path of the form
$(a_1, h, h', v)$, where $h, h'\in V(A^*)$ and $v\in V$, can act as a
witness for at most $2\ell-1$ bad $c_r$'s, since $v$ is in the vertex
set of at most $2\ell-1$ copies of $H_1^{\ell-1}[2\ell-1]$; this
implies that $v$ is in at most $2\ell-1$ $V(T_r)$'s. Starting at
$a_1$, for a blue-red-green path (resp.~green-red-blue path), there
are at most $\Delta_2$ (resp.~$\Delta_3$) choices for $h$, $\ell-1$
choices for $h'$, and $\Delta_3$ (resp.~$\Delta_2$) choices for $v$;
hence, there are at most $2\ell\Delta_2\Delta_3$ witnesses, and
$2\ell(2\ell-1)\Delta_2\Delta_3$ is an upper bound on the number of
$c_r$'s that are bad via this case.

\item[(2)]

Case 2: $\hat{H}_1^{\ell-1}$ forms with $G_2$ an alternating-$(\ell,
2, \ell, 2)$-bag via $D = D_{1, l^*}$ for some $l^*\in [\ell]$ and
some copy $A^*$ of $K_\ell(G_1)$ disjoint from $D$.

If $D$ forms an alternating-$(\ell, 2, \ell, 2)$-bag with $A^*$, then
$C_1$ forms an alternating-$(2\ell-1, 2, \ell, 2)$-bag with $A^*$,
which implies the existence of either a blue-red-green path or a
green-red-blue path $(w, h, h', c)$, where $w\in V(C_1)$ and $h,
h'\in{V(A^*)}$ that we think of as a witness to $c$ being bad. Each
blue-red-green path and green-red-blue path of the form $(w, h, h',
v)$, where $w\in V(C_1)$, $h, h'\in V(A^*)$ and $v\in V$, can act as a
witness for at most one bad $c_r$ (exactly one if $v = c_r$ for some
$k$). For a blue-red-green path (resp.~green-red-blue path), there are
at most $2\ell-1$ choices for $w$, $\Delta_2$ (resp.~$\Delta_3$)
choices for $h$, $\ell-1$ choices for $h'$, and $\Delta_3$
(resp.~$\Delta_2$) choices for $v$; hence, there are at most
$2\ell(2\ell-1)\Delta_2\Delta_3$ witnesses, and
$2\ell(2\ell-1)\Delta_2\Delta_3$ is an upper bound on the number of
$c_r$'s that are bad via this case.

\end{itemize}

Suppose we add $(a_1, a_2)$ to $E_3$.

\begin{itemize}

\item[(3)]

Case 3: $\hat{H}_1^{\ell-1}$ forms with $G_2$ an alternating-$(\ell,
2, \ell, 2)$-bag via $R = R_{r^*, l_1}$ for some $l_1\in [\ell]$ and
$D = D_{2, l_2}$ for some $l_2\in [\ell]$.

If $R$ forms an alternating-$(\ell, 2, \ell, 2)$-bag with $D$, then
$T_{r^*}$ forms an alternating-$(2\ell-1, 2, 2\ell-1, 2)$-bag with
$C_2$, which implies the existence of a blue-red-green path $(a_1,
a_2, w, u)$, where $w\in V(C_2)$ and $u\in V(T_{r^*})$ that we think
of as a witness to $c$ being bad. Each blue-red-green path of the form
$(a_1, a_2, w, v)$, where $w\in V(C_2)$ and $v\in V$, can act as a
witness for at most $2\ell-1$ bad $c_r$'s, since $v$ is in the vertex
set of at most $2\ell-1$ copies of $H_1^{\ell-1}[2\ell-1]$; this
implies that $v$ is in at most $2\ell-1$ $V(T_r)$'s. There are at most
$2\ell-1$ choices for $w$ and at most $\Delta_3$ choices for $v$;
hence, there are at most $(2\ell-1)\Delta_3$ witnesses, and
$(2\ell-1)^2\Delta_3$ is an upper bound on the number of $c_r$'s that
are bad via this case.

\end{itemize}

Summing over the three cases, we obtain that the total number of bad
$c_r$'s is at most
$4\ell(2\ell-1)\Delta_2\Delta_3+(2\ell-1)^2\Delta_3$, which is at most
$5\ell(2\ell-1)\Delta_2\Delta_3$ (assuming
$\Delta_2 > \lfloor (2\ell-1)/\ell\rfloor = 1$). Therefore, to show
that there exists an $\hat{H}_1^{\ell-1}$ that does not form an
alternating-$(\ell, 2, \ell, 2)$-bag with $G_2$, it suffices to show
\[5\ell(2\ell-1)\Delta_2\Delta_3+1\leq q,
\]
which is assumed. Choose $H_2^{\ell-1}$ to be such an
$\hat{H}_1^{\ell-1}$, and we are done.
\end{proof}

We now prove Theorem~\ref{thm:problem-1}.

\begin{proof}[Proof of Theorem~\ref{thm:problem-1}]
By Lemma~\ref{lemma:ham-cycle}, there exists a copy $H_0^{\ell-1}$ of
$H^{\ell-1}$ in $G_1$. To show that there exists a copy of
$H^{\ell-1}$ in $G_1$ that does not form an alternating-$(\ell, 2,
\ell, 2)$-bag with $G_2$ and $G_3$, we proceed by induction. Let
$G_{2, 0} = G_2$ and $G_{3, 0} = G_3$ as given. If $E_2 = E_{2, 0}$
and $E_3 = \varnothing$, then $H^{\ell-1}_0$ does not form an
alternating-$(\ell, 2, \ell, 2)$-bag with $G_2$ and $G_3$. For the
inductive step, let us add edges to $E_3$ until $H^{\ell-1}_0$ forms
an alternating-$(\ell, 2, \ell, 2)$-bag with $G_2$ and $G_3$. Let $e =
\edge{a_1}{a_2}$ be the last edge added to $E_3$, and let $G_3' = (V,
E_3\backslash\{e\})$ be $G_3$ just prior to adding $e$, so that
$H^{\ell-1}_0$ does not form an alternating-$(\ell, 2, \ell, 2)$-bag
with $G_2$ and $G_3'$. Then, by Lemma~\ref{lemma:no-bag}, there exists
a copy of $H^{\ell-1}$ in $G_1$ that does not form an
alternating-$(\ell, 2, \ell, 2)$-bag with $G_2$ and $G_3$. This
contradiction indicates that there exists a copy of $H^{\ell-1}$ in
$G_1$ that does not form an alternating-$(\ell, 2, \ell, 2)$-bag with
$G_{2, 0}$ and $G_{3, 0}$.
  
\end{proof}

Finally, we use Theorem~\ref{thm:problem-1} to prove
Theorem~\ref{thm:baranyai-1}, which we restate as follows.

\vspace{6pt}
\noindent{\bf Theorem~\ref{thm:baranyai-1}.}~
{\em
Let $\alpha,\beta\in(0,1)$ such that $\alpha+\beta\geq 1$.  
Let $n$ be a positive integer. Let $k$ be $O(1)$ and $\ell$ be
$o(\sqrt{n})$. There exists a cyclic ordering of the $k$-element
subsets of $[n]$ with the following properties:

\begin{itemize}
\item[(1)]
any $\ell$ consecutive $k$-element subsets in the ordering are
disjoint, and in particular, form a $(k, \ell)$-parpartition of $[n]$,
which we refer to as a consecutive-$(k, \ell)$-parpartition; and
\item[(2)]
given any consecutive-$(k, \ell)$-parpartition $P_1$, no disjoint
consecutive-$(k, \ell)$-parpartition $P_2$ is too close to $P_1$
(that is, $P_1$ and $P_2$ are not $(\alpha,\beta)$-close).
\end{itemize}
}

\begin{proof}[Proof of Theorem~\ref{thm:baranyai-1}]
Let $V$ be a set of ${n\choose k}$ vertices corresponding bijectively
to the ${n\choose k}$ $k$-element subsets of $[n]$. Let $E_1$
be the set of edges $\edge{u}{v}$, where $u, v\in V$,
such that the $k$-element subsets corresponding to $u$ and $v$ are
disjoint. Let $E_2$
(resp.~$E_3$) be the set of edges $\edge{u}{v}$, where $u, v\in V$,
such that the $k$-element subsets corresponding to $u$ and $v$ 
intersect in more than $\alpha{k}$ (resp.~$\beta{k}$)
elements. Letting $G_1 = (V, E_1)$, $G_2 = (V, E_2)$, and $G_3 = (V,
E_3)$, we show that $G_1$, $G_2$, and $G_3$ satisfy the two
conditions in Theorem~\ref{thm:problem-1} with $q = m$.

Recall that Condition~\ref{equation:condition-1} from Theorem~\ref{thm:problem-1} states that.

\[m\left(\frac{(2\ell-1)^2-1}{(2\ell-1)^2}\right)+\frac{4\ell-3+q}{(2\ell-1)^2}\leq\delta_1,
\]
where $q\geq 1$. Let us choose $q = m$. Then we can rewrite \ref{equation:condition-1} as
\[m\left(\frac{(2\ell-1)^2-1}{(2\ell-1)^2}+\eps\right)\leq\delta_1,
\]
where $0\leq\eps\leq 1/(2\ell-1)^2$, assuming $m >
(4\ell-3)/((2\ell-1)^2\eps-1)$ holds. Let us choose $\eps =
1/2(2\ell-1)^2$. We have $m = {n\choose k}$ and $\delta_1 =
{n-k\choose k}$. Then this becomes
\[{n\choose k}\left(1-\frac{1}{2(2\ell-1)^2}\right)\leq {n-k\choose k}.
\]
As $n$ grows large, ${n-k\choose k}/{n\choose k}$ approaches $1$, and
so the inequality holds for adequately large $n$.

Now recall Condition~\ref{equation:condition-2}, which states that
\[\Delta_2\Delta_3\leq\frac{m-1}{5\ell(2\ell-1)}.
\]
We have that the degree of every vertex in $V(G_2)$ is
\[\sum_{1\leq i < (1-\alpha)k}{k\choose i}{n-k\choose i}, 
\]
which is $O(n^{\lfloor(1-\alpha)k\rfloor})$
or $O(n^{(1-\alpha)k-1})$ if $(1-\alpha)k$ is an integer.
Similarly, the degree of every vertex in
$V(G_3)$ is
\[\sum_{1\leq i < (1-\beta)k}{k\choose i}{n-k\choose i},
\]
which is $O(n^{\lfloor(1-\beta)k\rfloor})$ or $O(n^{(1-\beta)k-1})$ if
$(1-\beta)k$ is an integer.  So $\Delta_2\Delta_3$ is $O(n^{k-1})$
whenever $\alpha+\beta\geq{1}$.  In this case, the right side of \ref{equation:condition-2}
is $\omega(n^{k-1})$ (as $\ell$ is $o(\sqrt{n})$), so the inequality
holds asymptotically.

Thus, by Theorem~\ref{thm:problem-1}, for some Hamiltonian cycle $H$
there exists $H^{\ell-1}$ in $G_1$ such that no two disjoint copies of
$K_\ell$ in $H^{\ell-1}$ form an alternating-$(\ell, 2, \ell, 2)$-bag
with $G_2$ and $G_3$. Observe that (1) the copies of $K_\ell$ in
$H^{\ell-1}$ are exactly the subgraphs of $H^{\ell-1}$ induced by
$\ell$ consecutive vertices on $H^{\ell-1}$, which correspond to a
$(k, \ell)$-parpartition, and (2) an alternating-$(\ell, 2, \ell,
2)$-bag between two disjoint copies of $K_\ell$ correspond to the two
corresponding $(k, \ell)$-parpartitions being $(\alpha,
\beta)$-close. Therefore, $H$ is an ordering of the $k$-element
subsets of $[n]$ that has the two properties in
Theorem~\ref{thm:baranyai-1}, and we are done.

\end{proof}


\section{Conclusion}
In this paper, we have studied variants of Baranyai's theorem under
an additional restriction. In particular, we have shown that, given
$k$, $\ell$, and $n$ satisfying $k^2\ell\leq n/3$, one can find
$\lfloor{n\choose k}/\ell\rfloor$ ways to partition the integer set
$[n]$ into $k$-element subsets such that no two parpartitions have
significant overlap. This greatly increases the applicability of a
recent parallel result by Katona and Katona
(Theorem~\ref{thm:katona-katona-1}) that only handles $k$, $\ell$,
and $n$ such that $k^2\ell = o(\sqrt{n})$.

We have also given a result on conditions sufficient for a graph to
contain a Hamiltonian cycle to the $p$th power of a certain property,
which answers a question posed by Katona and Katona, and used this to
prove another result that can also be viewed as improving on Katona's
and Katona's Theorem \ref{thm:katona-katona-1}.

We pose here some directions for future research.

First, we are interested in extending Theorem
\ref{thm:modified-katona} to parpartitions of larger size.

\begin{problem}
Given positive integers $k$, $\ell$, and $n$ such that $k\ell = n$, can one find
${n\choose k}/\ell$ parpartitions of $[n]$ such that no $k$-element
set appears in two parpartitions and no two parpartitions are too
close?
\end{problem}

It would also be interesting to study Baranyai's theorem with a different condition on the relationship between partitions (or parpartitions).

\begin{problem}
Given positive integers $k$ and $n$ such that $k$ divides $n$, find a relation between partitions of $[n]$ into $k$-element subsets for which we can decompose the family of $k$-element sets of $[n]$ into partitions of $[n]$ such that no subset of these partitions satisfy this relation.
\end{problem}

For instance, it may be interesting to consider a different definition of closeness of partitions, e.g., involving more than two partitions or more than two pairs of sets in the partitions.

\section{Acknowledgements}

This research was conducted at the University of Minnesota Duluth REU with support from Jane Street Capital, NSF Grant 1949884, and donations from Ray Sidney and Eric Wepsic. The author thanks Colin Defant for suggesting the problem and Maya Sankar for extensive feedback during the editing process. The author also thanks Joe Gallian, Colin Defant, Noah Kravitz, Mitchell Lee, and Maya Sankar for their guidance and support.

\end{document}